\documentclass[12pt,reqno]{amsart}
\usepackage{amssymb}
\usepackage{amsxtra}
\usepackage[mathscr]{eucal}

\newcommand{\Z}{{\mathbf Z}}

\newcommand{\Q}{{\mathbf Q}}

\newcommand{\eop}{\hfill$\square$}

\theoremstyle{plain}
\newtheorem{Thm}{Theorem}

\newtheorem{Cor}{Corollary}
\newtheorem{Lem}{Lemma}
\newtheorem{Prop}{Proposition}
\theoremstyle{definition}

\theoremstyle{remark}

\begin{document}

\title[Signed $q$-Analogs of Tornheim's Double Series]{Signed $q$-Analogs of Tornheim's Double Series}

\date{\today}

\author{Xia Zhou}
\address{Department of Mathematics\\ Zhejiang University, Hangzhou,
310027, P.\ R.\ China} \email[]{xiazhou0821@hotmail.com}

\author{Tianxin Cai}
\address{Department of Mathematics, Zhejiang University, Hangzhou,
310027, P.\ R.\ China} \email[]{txcai@mail.hz.zj.cn}

\author{David~M. Bradley}
\address{Department of Mathematics \& Statistics\\
         University of Maine\\
         5752 Neville Hall
         Orono, Maine 04469-5752\\
         U.S.A.}
\email[]{bradley@math.umaine.edu, dbradley@member.ams.org}



\begin{abstract}
We introduce signed $q$-analogs of Tornheim's double series, and
evaluate them in terms of double $q$-Euler sums.  As a consequence,
we provide explicit evaluations of signed and unsigned Tornheim
double series, and correct some mistakes in the literature.
\end{abstract}

\maketitle

%\tableofcontents
\interdisplaylinepenalty=500

\section{Introduction}\label{sect:Intro}
 Let $k$ be a positive integer.  Sums of the form
\begin{equation}\label{MzvDef}
   \zeta(s_1,s_2,\dots,s_k) := \sum_{n_1>n_2>\cdots>n_k>0}\;
   \prod_{j=1}^k n_j^{-s_j},
   \quad \sum_{j=1}^m \Re(s_j) > m,\quad m=1,2,\dots,k.
\end{equation}
have attracted increasing attention in recent years; see
eg.~\cite{ParmEuler,
GoldVar,BBB,BBBLc,BBBLa,BowBrad3,BowBrad1,BowBradRyoo,Prtn,DBqKarl,BK1,Hoff02,IT,LeM,Tera1,Tera2,Tera3}.
The survey
articles~\cite{BowBradSurvey,Cartier,Cartier06,Wald02,Wald04,Zud}
provide an extensive list of references.  In~\eqref{MzvDef} the sum
is over all positive integers $n_1,n_2,\dots,n_k$ satisfying the
indicated inequalities.  Of course~\eqref{MzvDef} reduces to the
familiar Riemann zeta function when $k=1$.   When the arguments are
all positive integers, we refer to~\eqref{MzvDef} as a multiple zeta
value and note that in this case, $s_1>1$ is necessary and
sufficient for convergence.

The problem of evaluating $\zeta(s_1,s_2)$ with integers $s_1>1$ and
$s_2>0$ seems to have been first proposed in a letter from Goldbach
to Euler~\cite{LE2} in 1742. (See also~\cite{LE,Goldbach}
and~\cite[p.\ 253]{Berndt1}.)  Calculating several examples led
Euler to infer a closed form evaluation in terms of values of the
Riemann zeta function in the case when $s_1+s_2$ is odd.
In~\cite{BBB}, Borwein, Bradley and Broadhurst considered the more
general Euler sum
\begin{equation}\label{EulerSum}
   \zeta(s_1,s_2,\dots,s_k;\sigma_1,\sigma_2,\dots,\sigma_k) :=
   \sum_{n_1>n_2>\cdots>n_k>0}\; \prod_{j=1}^k \sigma_j^{n_j} n_j^{-s_j}
\end{equation}
with each $\sigma_j\in\{-1,1\}$.   Among the many other results
for~\eqref{EulerSum} listed therein is an explicit formula for the
case $k=2$ that reduces to Euler's evaluation when
$\sigma_1=\sigma_2=1$. For the case of arbitrary $k$, several
infinite classes of closed form evaluations for~\eqref{MzvDef}
and~\eqref{EulerSum} are proved
in~\cite{BBBLc,BBBLa,BowBrad3,BowBrad1,BowBradRyoo}. Subsequently,
Bradley~\cite{DBqMzv,DBqKarl,DBqSum,DBqDecomp} found $q$-analogs
of~\eqref{MzvDef} and~\eqref{EulerSum} and thoroughly investigated
their properties.

In~\cite{Torn} Tornheim introduced the double series
\[
   T(r,s,t) := \sum_{u,v=1}^\infty \frac{1}{u^{r} v^{s} (u+v)^{t}},
\]
with nonnegative integers $r$, $s$ and $t$ satisfying $r+t>1$,
$s+t>1$ and $r+s+t>2$. Huard, Williams and Zhang~\cite{HWZ}
evaluated Torheim's series in terms of sums of products of Riemann
zeta values when $r+s+t$ is odd. Subbarao and
Sitaramachandrarao~\cite{SitSub} considered the alternating variants
\[
   R(r,s,t) := \sum_{u,v=1}^\infty \frac{(-1)^v}{u^rv^s(u+v)^t},
   \qquad
    S(r,s,t) := \sum_{u,v=1}^\infty
   \frac{(-1)^{u+v}}{u^rv^s(u+v)^t},
\]
and posed the problem to evaluate $R(r,r,r)$ and $S(r,r,r)$ for any
positive integer $r$.  Tsumura~\cite[Cor.~3]{Tsu1} and~\cite[Theorem
3.6]{Tsu2} tackled this problem and later~\cite{Tsu3} evaluated
$S(r,s,t)$ for any positive integers $r$, $s$, $t$ such that $r+s+t$
is odd. Tsumura's method is elementary but complicated, and also has
some mistakes. In this paper, we give a simple formula that
expresses $R$, $S$ and $T$ in terms of double Euler sums. In light
of the aforementioned formula of Borwein, Bradley and Broadhurst for
the double Euler sums, our results yield a closed form evaluation
for the series $R$, $S$ and $T$ whenever the arguments $r$, $s$ and
$t$ are positive integers with $r+s+t$ odd. More generally, we
consider $q$-analogs of $R$, $S$, $T$ and show how they may all be
evaluated in terms of double $q$-Euler sums.

\section{q-Analogs}
Henceforth assume $q$ is real and $q>1$. The $q$-analog of a
positive integer $n$ is
\[
   [n]_q := \sum_{j=0}^{n-1} q^j = \frac{q^n-1}{q-1}.
          %= \begin{cases}\displaystyle\frac{q^n-1}{q-1} &\mbox{if $q>1$},\\
           %              n &\mbox{if $q=1$}.
           % \end{cases}
\]
Let $k$ be a positive integer, let $s_1,s_2,\dots,s_k$ be real
numbers, and let $\sigma_1,\sigma_2,\dots,\sigma_k\in\{-1,1\}$.
Define the $q$-Euler sum
\begin{equation}\label{qEulerSum}
   \zeta_q[s_1,s_2,\dots,s_k;\sigma_1,\sigma_2,\dots,\sigma_k] :=
   \sum_{n_1>n_2>\cdots>n_k>0}\;\prod_{j=1}^k
   \frac{\sigma_j^{n_j}q^{(s_j-1)n_j}}{[n_j]_q^{s_j}}
\end{equation}
and note that this coincides with the special case
$\mathrm{Li}_{s_1,\dots,s_k}[\sigma_1 q^{s_1-1},\dots,\sigma_k
q^{s_k-1}]$ of the multiple
$q$-polylogarithm~\cite[eq.~(6.2)]{DBqMzv}.  If $\sigma_j=1$ for
each $j=1,2,\dots,k$, then we recover the multiple $q$-zeta value
$\zeta[s_1,s_2,\dots,s_k]$ of~\cite{DBqMzv,DBqSum,DBqDecomp}.

Let $\sigma,\tau\in\{-1,1\}$.  Define $q$-analogs of the signed and
unsigned Tornheim double series $R$, $S$ and $T$ by
\[
   T[r,s,t;\sigma,\tau] := \sum_{u,v=1}^\infty
   \frac{\sigma^u\tau^v
   q^{(r+t-1)u+(s+t-1)v}}{[u]_q^r[v]_q^s[u+v]_q^t}.
\]
The sum
\[
   \varphi[s;\sigma] := \sum_{n=1}^\infty \frac{(n-1)\sigma^nq^{(s-1)n}}{[n]_q^s}
   = \sum_{n=1}^\infty \frac{n \sigma^n q^{(s-1)n}}{[n]_q^s} -
   \zeta_q[s;\sigma],
\]
the case $\sigma=1$ of which was defined in~\cite{DBqDecomp}, will
also be needed.   As in~\cite{BBB}, it is convenient to combine
signs and exponents in~\eqref{EulerSum} and~\eqref{qEulerSum} into a
single list by writing $s_j$ if $\sigma_j=1$ and $\overline{s_j}$ if
$\sigma_j=-1$. For consistency one may do this also for $T$ and
$\varphi$; thus for example, $\varphi[s;1]=\varphi[s]$,
$\varphi[s;-1]=\varphi[\,\overline{s}\,]$,
$R(r,s,t)=T(r,\overline{s},t)$ and
$S(r,s,t)=T(\overline{r},\overline{s},t)$.  We also employ the
notation
\[
   \binom{z}{a,\, b} := \binom{z}{a}\binom{z-a}{b} =
   \binom{z}{b}\binom{z-b}{a}
\]
for the trinomial coefficient, in which $a$, $b$ are nonnegative
integers, and which reduces to $z!/a!b!(z-a-b)!$ if $z$ is also an
integer exceeding $a+b$.  

The following theorem shows how the $q$-analogs of $R$, $S$ and $T$
are related to the $q$-Euler sums.

\begin{Thm}\label{thm:qT} Let $r$ and $s$ be positive integers, and let $t$ be a real number. Then
  \begin{align*}
     T[r,s,t] &=
     \sum_{a=0}^{r-1}\sum_{b=0}^{r-1-a}\binom{a+s-1}{a,\, b}(1-q)^b\zeta_q[s+t+a,r-a-b]\\
     &+\sum_{a=0}^{s-1}\sum_{b=0}^{s-1-a}\binom{a+r-1}{a,\, b}(1-q)^b\zeta_q[r+t+a,s-a-b]\\
     &-\sum_{j=1}^{\mathrm{min}(r,s)}\binom{r+s-j-1}{r-j,\, s-j}(1-q)^j\varphi[r+s+t-j],\\
    T[\overline{r},\overline{s},t]
    &=\sum_{a=0}^{r-1}\sum_{b=0}^{r-1-a}\binom{a+s-1}{a,\, b}(1-q)^b\zeta_q[\,\overline{s+t+a},r-a-b]\\
    &+\sum_{a=0}^{s-1}\sum_{b=0}^{s-1-a}\binom{a+r-1}{a,\, b}(1-q)^b\zeta_q[\,\overline{r+t+a},s-a-b]\\
    &-\sum_{j=1}^{\mathrm{min}(r,s)}\binom{r+s-j-1}{r-j,\, s-j}(1-q)^j\varphi[\,\overline{r+s+t-j}\,],\\
   T[r,\overline{s},t]
   &=\sum_{a=0}^{r-1}\sum_{b=0}^{r-1-a}\binom{a+s-1}{a,\, b}(1-q)^b
                     \zeta_q[\,\overline{s+t+a},\overline{r-a-b}\,]\\
   &+\sum_{a=0}^{s-1}\sum_{b=0}^{s-1-a}\binom{a+r-1}{a,\, b}(1-q)^b\zeta_q[\,r+t+a,\overline{s-a-b}\,]\\
   &-\sum_{j=1}^{\mathrm{min}(r,s)}\binom{r+s-j-1}{r-j,\, s-j}(1-q)^j(1+q)^{j-r-s-t}\zeta_{q^2}[r+s+t-j].
  \end{align*}
\end{Thm}

Taking the limit as $q\to 1+$ in Theorem~\ref{thm:qT} and noting the
restrictions on $r$, $s$ and $t$ now needed for convergence yields
the following
\begin{Cor}\label{cor:Torn2Euler} Let $r$ and $s$ be positive integers and let $t$ be a real number.
If $r+t>1$ and $s+t>1$, then
\begin{align*}
  T(r,s,t) &=
  \sum_{a=0}^{r-1}\binom{a+s-1}{s-1}\zeta(s+t+a,r-a)+\sum_{a=0}^{s-1}\binom{a+r-1}{r-1}\zeta(r+t+a,s-a);\\
  \intertext{if $r+t>0$ and $s+t>0$, then}
  S(r,s,t) &=\sum_{a=0}^{r-1}\binom{a+s-1}{s-1}\zeta(\overline{s+t+a},r-a)
             +\sum_{a=0}^{s-1}\binom{a+r-1}{r-1}\zeta(\overline{r+t+a},s-a);\\
  \intertext{if $r+t>1$ and $s+t>0$, then}
  R(r,s,t) &= \sum_{a=0}^{r-1}\binom{a+s-1}{s-1}\zeta(\overline{s+t+a},\overline{r-a})
             +\sum_{a=0}^{s-1}\binom{a+r-1}{r-1}\zeta(r+t+a,\overline{s-a}).
\end{align*}
\end{Cor}

Putting $t=0$ in Theorem~\ref{thm:qT} yields the following
decomposition formulas, the first of which was given under slightly
more restrictive hypotheses in~\cite[Theorem 2.1]{DBqDecomp}.
\begin{Cor} If $r$ and $s$ are positive integers, then
\begin{align*}
  \zeta_q[r]\zeta_q[s] &= \sum_{a=0}^{r-1}\sum_{b=0}^{r-1-a}\binom{a+s-1}{s-1}\binom{s-1}{b}(1-q)^b\zeta_q[s+a,r-a-b]\\
  &+\sum_{a=0}^{s-1}\sum_{b=0}^{s-1-a}\binom{a+r-1}{r-1}\binom{r-1}{b}(1-q)^b\zeta_q[r+a,s-a-b]\\
  &-\sum_{j=1}^{\mathrm{min}(r,s)}\binom{r+s-j-1}{r-j,\, s-j}(1-q)^j\varphi[r+s-j];\\
  \zeta_q[\,\overline{r}\,]\zeta_q[\,\overline{s}\,]
  &=\sum_{a=0}^{r-1}\sum_{b=0}^{r-1-a}\binom{a+s-1}{s-1}\binom{s-1}{b}(1-q)^b\zeta_q[\,\overline{s+a},r-a-b]\\
  &+\sum_{a=0}^{s-1}\sum_{b=0}^{s-1-a}\binom{a+r-1}{r-1}\binom{r-1}{b}(1-q)^b\zeta_q[\,\overline{r+a},s-a-b]\\
  &-\sum_{j=1}^{\mathrm{min}(r,s)}\binom{r+s-j-1}{r-j,\, s-j}(1-q)^j\varphi[\,\overline{r+s-j}\,];\\
  \zeta_q[r]\zeta_q[\,\overline{s}\,]
  &= \sum_{a=0}^{r-1}\sum_{b=0}^{r-1-a}\binom{a+s-1}{s-1}\binom{s-1}{b}(1-q)^b\zeta_q[\,\overline{s+a},\overline{r-a-b}\,]\\
  &+\sum_{a=0}^{s-1}\sum_{b=0}^{s-1-a}\binom{a+r-1}{r-1}\binom{r-1}{b}(1-q)^b\zeta_q[r+a,\overline{s-a-b}\,]\\
  &-\sum_{j=1}^{\mathrm{min}(r,s)}\binom{r+s-j-1}{r-j,\, s-j}(1-q)^j(1+q)^{j-r-s}
               \zeta_{q^2}[r+s-j].
\end{align*}
\end{Cor}

Taking the limit as $q\to 1+$ in Corollary 2 and noting the
additional restrictions needed on $r$ and $s$ to guarantee
convergence in this case yields the following decomposition
formulas, the first of which was known to Euler.
\begin{Cor}
   If $r-1$ and $s-1$ are positive integers, then
  \begin{align*}
    \zeta(r)\zeta(s) &=
    \sum_{a=0}^{r-1}\binom{a+s-1}{s-1}\zeta(s+a,r-a)
     +\sum_{a=0}^{s-1}\binom{a+r-1}{r-1}\zeta(r+a,s-a);\\
     \intertext{if $r$ and $s$ are positive integers, then}
    \zeta(\overline{r})\zeta(\overline{s})
    &=\sum_{a=0}^{r-1}\binom{a+s-1}{s-1}\zeta(\overline{s+a},r-a)
    +\sum_{a=0}^{s-1}\binom{a+r-1}{r-1}\zeta(\overline{r+a},s-a);\\
   \intertext{if $r-1$ and $s$ are positive integers, then}
   \zeta(r)\zeta(\overline{s})
   &= \sum_{a=0}^{r-1}\binom{a+s-1}{s-1}\zeta(\overline{s+a},\overline{r-a})
     +\sum_{a=0}^{s-1}\binom{a+r-1}{r-1}\zeta(r+a,\overline{s-a}).
  \end{align*}
\end{Cor}

\section{Proof of Theorem~\ref{thm:qT}.}
The key ingredient is the following partial fraction decomposition.

\begin{Lem}\label{lem:qParfrac}
  If $r$, $s$, $u$, and $v$ are all positive integers, then
  \begin{align*}
     \frac{1}{[u]_q^r[v]_q^s}
     &= \sum_{a=0}^{r-1}\sum_{b=0}^{r-1-a}\binom{a+s-1}{a,\, b}
         \frac{(1-q)^b\, q^{(s-1-b)u+av}}{[u]_q^{r-a-b}[u+v]_q^{s+a}}\\
     &+\sum_{a=0}^{s-1}\sum_{b=0}^{s-1-a}\binom{a+r-1}{a,\, b}
       \frac{(1-q)^b\, q^{au+(r-1-b)v}}{[v]_q^{s-a-b}[u+v]_q^{r+a}}\\
     &-\sum_{j=1}^{\min(r,s)}\binom{r+s-j-1}{r-j,\, s-j}\frac{(1-q)^j\, q^{(s-j)u+(r-j)v}}{[u+v]_q^{r+s-j}}.
  \end{align*}
\end{Lem}

\noindent{\bf Proof.}  Let $x$ and $y$ be nonzero real numbers such
that $x+y+(q-1)xy \ne 0$.  As in~\cite{DBqDecomp}, observe that if
we apply the partial differential operator
\[
   \frac{1}{(r-1)!}\bigg(-\frac{\partial}{\partial x}\bigg)^{r-1}
   \frac{1}{(s-1)!}\bigg(-\frac{\partial}{\partial y}\bigg)^{s-1}
\]
to both sides of the identity
\[
   \frac{1}{xy}
 = \frac{1}{x+y+(q-1)xy}\bigg(\frac{1}{x}+\frac{1}{y}+q-1\bigg),
\]
then we obtain the identity~\cite[Lemma 3.1]{DBqDecomp}
\begin{align*}
   \frac{1}{x^r y^s}
   &= \sum_{a=0}^{r-1}\,\sum_{b=0}^{r-1-a}
      \binom{a+s-1}{a,\, b}
      \frac{(1-q)^b(1+(q-1)y)^a(1+(q-1)x)^{s-1-b}}
           {x^{r-a-b}(x+y+(q-1)xy)^{s+a}}\\
   &+ \sum_{a=0}^{s-1}\,\sum_{b=0}^{s-1-a}
      \binom{a+r-1}{a,\, b}
      \frac{(1-q)^b(1+(q-1)x)^a(1+(q-1)y)^{r-1-b}}
           {y^{s-a-b}(x+y+(q-1)xy)^{r+a}}\\
   &- \sum_{j=1}^{\min(r,s)}\binom{r+s-j-1}{r-j,\, s-j}
      \frac{(1-q)^j\, (1+(q-1)y)^{r-j}(1+(q-1)x)^{s-j}}{(x+y+(q-1)xy)^{r+s-j}}.
\end{align*}
Now let $x=[u]_q$, $y=[v]_q$ and note that then $1+(q-1)x=q^u$,
$1+(q-1)y=q^v$ and $x+y+(q-1)xy=[u+v]_q$. \eop

To prove Theorem~\ref{thm:qT}, multiply both sides of
Lemma~\ref{lem:qParfrac} by
\[
   \frac{\sigma^u\tau^v q^{(r+t-1)u+(s+t-1)v}}{[u+v]_q^t}
\]
and sum over all ordered pairs of positive integers $(u,v)$ to
obtain
\begin{align*}
  T[r,s,t;\sigma,\tau] &=
  \sum_{a=0}^{r-1}\sum_{b=0}^{r-1-a}\binom{a+s-1}{a,\, b}(1-q)^b
  \sum_{u,v=1}^\infty\frac{\sigma^u\tau^vq^{(r-a-b-1)u}q^{(s+t+a-1)(u+v)}}{[u]_q^{r-a-b}[u+v]_q^{s+t+a}}\\
  &+\sum_{a=0}^{s-1}\sum_{b=0}^{s-1-a}\binom{a+r-1}{a,\ b}(1-q)^b
  \sum_{u,v=1}^\infty\frac{\sigma^u\tau^vq^{(s-a-b-1)v}q^{(r+t+a-1)(u+v)}}{[v]_q^{s-a-b}[u+v]_q^{r+t+a}}\\
  &-\sum_{j=1}^{\min(r,s)}\binom{r+s-j-1}{r-j,\, s-j}(1-q)^j
  \sum_{u,v=1}^\infty \frac{\sigma^u\tau^v q^{(r+s+t-j-1)(u+v)}}{[u+v]_q^{r+s+t-j}}.
\end{align*}
It follows that
\begin{align*}
  T[r,s,t;\sigma,\sigma] &=\sum_{a=0}^{r-1}\sum_{b=0}^{r-1-a}\binom{a+s-1}{a,\, b}(1-q)^b
  \sum_{m>n>0}\;\frac{\sigma^mq^{(s+t+a-1)m}q^{(r-a-b-1)n}}{[m]_q^{s+t+a}[n]_q^{r-a-b}}\\
  &+\sum_{a=0}^{s-1}\sum_{b=0}^{s-1-a}\binom{a+r-1}{a,\, b}(1-q)^b
  \sum_{m>n>0}\;\frac{\sigma^mq^{(r+t+a-1)m}q^{(s-a-b-1)n}}{[m]_q^{r+t+a}[n]_q^{s-a-b}}\\
  &-\sum_{j=1}^{\min(r,s)}\binom{r+s-j-1}{r-j,\, s-j}(1-q)^j
  \sum_{m>n>0}\;\frac{\sigma^mq^{(r+s+t-j-1)m}}{[m]_q^{r+s+t-j}}\\
  &=\sum_{a=0}^{r-1}\sum_{b=0}^{r-1-a}\binom{a+s-1}{a,\, b}(1-q)^b\zeta_q[s+t+a,r-a-b;\sigma,1]\\
  &+\sum_{a=0}^{s-1}\sum_{b=0}^{s-1-a}\binom{a+r-1}{a,\, b}(1-q)^b\zeta_q[r+t+a,s-a-b;\sigma,1]\\
  &-\sum_{j=1}^{\min(r,s)}\binom{r+s-j-1}{r-j,\, s-j}(1-q)^j\varphi[r+s+t-j;\sigma],
\end{align*}
and also that
\begin{align*}
  T[r,s,t;1,-1]
  &=\sum_{a=0}^{r-1}\sum_{b=0}^{r-1-a}\binom{a+s-1}{a,\, b}(1-q)^b
  \sum_{m>n>0}\;\frac{(-1)^mq^{(s+t+a-1)m}(-1)^nq^{(r-a-b-1)n}}{[m]_q^{s+t+a}[n]_q^{r-a-b}}\\
  &+\sum_{a=0}^{s-1}\sum_{b=0}^{s-1-a}\binom{a+r-1}{a,\, b}(1-q)^b
  \sum_{m>n>0}\;\frac{q^{(r+t+a-1)m}(-1)^nq^{(s-a-b-1)n}}{[m]_q^{r+t+a}[n]_q^{s-a-b}}\\
  &-\sum_{j=1}^{\min(r,s)}\binom{r+s-j-1}{r-j,\, s-j}(1-q)^j
  \sum_{m>n>0}\;\frac{(-1)^nq^{(r+s+t-j-1)m}}{[m]_q^{r+s+t-j}}\\
  &=\sum_{a=0}^{r-1}\sum_{b=0}^{r-1-a}\binom{a+s-1}{a,\, b}(1-q)^b\zeta_q[s+t+a,r-a-b;\sigma,1]\\
  &+\sum_{a=0}^{s-1}\sum_{b=0}^{s-1-a}\binom{a+r-1}{a,\, b}(1-q)^b\zeta_q[r+t+a,s-a-b;\sigma,1]\\
  &+\sum_{j=1}^{\min(r,s)}\binom{r+s-j-1}{r-j,\, s-j}(1-q)^j\sum_{m>0}\frac{q^{(r+s+t-j-1)2m}}{[2m]_q^{r+s+t-j}}.
\end{align*}
Since
\[
   [2m]_q = \frac{q^{2m}-1}{q-1} = \frac{q^2-1}{q-1}\cdot\frac{q^{2m}-1}{q^2-1}=(q+1)[m]_{q^2},
\]
the proof of Theorem~\ref{thm:qT} is complete. \eop

\section{Examples}

Again, let $\sigma,\tau\in\{-1,1\}$. It is known~\cite{BBB,FlajSalv}
that if $s$ and $t$ are positive integers such that $s+t$ is odd and
$s>(1+\sigma)/2$, then $\zeta(s,t;\sigma,\tau)$ lies in the
polynomial ring $\Q[\{\zeta(k;\pm 1): k\in\Z, 2\le k\le
s+t\}\cup\{\zeta(1;-1)\}]$. Since
\[
   \zeta(k;-1)= \zeta(\overline{k}) = \begin{cases} \big(2^{1-k}-1\big)\zeta(k) &\mbox{if
   $k>1$,} \\ -\log 2 &\mbox{if $k=1$,}\end{cases}
\]
it is clear that in fact, $\zeta(s,t;\sigma,\tau)\in \Q[\{\zeta(k):
k\in\Z, 2\le k\le s+t\}\cup\{-\log 2\}]$.   It follows that if
$r,s,t$ satisfy the conditions of Corollary~\ref{cor:Torn2Euler} and
if in addition $r+s+t$ is an odd integer, then the signed and
unsigned double Tornheim series $R(r,s,t)$, $S(r,s,t)$ and
$T(r,s,t)$ also lie in this ring.  It is possible to evaluate these
series explicitly if we recall the following formula from~\cite[eq.\
(75)]{BBB}.

\begin{Prop}\label{prop:BBB}
  Let $\sigma,\tau\in\{-1,1\}$, and
  let $s$ and $t$ be positive integers such that $s+t$ is odd,
  $s>(1+\sigma)/2$, and $t>(1+\tau)/2$.  Then
  \begin{multline}\label{BBBeq75}
    \zeta(s,t;\sigma,\tau)
    = \tfrac12\big(1+(-1)^s\big)\zeta(s;\sigma)\zeta(t;\tau)
    -\tfrac12\zeta(s+t;\sigma\tau)\\
    +(-1)^{t}\sum_{0\le k\le t/2}
    \binom{s+t-2k-1}{s-1}\zeta(2k;\sigma\tau)\zeta(s+t-2k;\sigma)\\
    +(-1)^{t}\sum_{0\le k\le s/2}\binom{s+t-2k-1}{t-1}\zeta(2k;\sigma\tau)\zeta(s+t-2k;\tau).
  \end{multline}
\end{Prop}
In Proposition 1, it is understood that $\zeta(0;\sigma\tau)=-1/2$
in accordance with the analytic continuation of $s\mapsto
\zeta(s;\sigma\tau)$.  The restriction $t>(1+\tau)/2$ can be removed
if in~\eqref{BBBeq75} we interpret $\zeta(1;1)=0$ wherever it
occurs.  That is, if $\sigma\in\{-1,1\}$ and $s$ is an even positive
integer, then
\begin{equation}\label{BBBeq75t=1}
   \zeta(s,1;\sigma,1) = \frac12
   (s-1)\zeta(s+1;\sigma)+\frac12\zeta(s+1)-\sum_{k=1}^{(s/2)-1}
   \zeta(2k;\sigma)\zeta(s+1-2k).
\end{equation}
The case $\sigma=1$ of~\eqref{BBBeq75t=1} is subsumed by another
formula~\cite[eq.\ (31)]{BBB} of Euler, namely
\begin{equation}\label{zs1}
   \zeta(s,1) = \frac{s}{2}\,\zeta(s+1) - \frac12 \sum_{k=2}^{s-1}
   \zeta(k)\zeta(s+1-k),
\end{equation}
which is valid for \emph{all} integers $s>1$, not just for even $s$.

In~\cite{Tsu3}, Tsumura listed evaluation formulas for $S(r,s,t)$
when $r+s+t\le 9$ is odd.  From Corollary~\ref{cor:Torn2Euler},
Proposition~\ref{prop:BBB} and equation~\eqref{BBBeq75t=1}, we can
deduce explicit formulas for $R(r,s,t)$, $S(r,s,t)$ and $T(r,s,t)$
when $r+s+t$ is odd.  In particular, we have the following new
results:
\begin{align*}
  R(1,1,1) &= -\frac58\zeta(3), & R(1,1,3) &=
  \frac1{16}\pi^2\zeta(3)-\frac{27}{32}\zeta(5),\\
  R(1,2,2) &= \frac{5}{48}\pi^2\zeta(3)-\frac32\zeta(5), & R(1,3,1)
  &= \frac1{12}\pi^2\zeta(3)-\frac{59}{32}\zeta(5),\\
  R(2,1,2) &= -\frac{5}{16}\pi^2\zeta(3)+\frac{107}{32}\zeta(5), &
  R(2,2,1) &= -\frac{5}{24}\pi^2\zeta(3)+\frac{59}{32}\zeta(5),\\
  R(3,1,1) &= \frac18\pi^2\zeta(3)-\frac{59}{32}\zeta(5),
\end{align*}
\begin{align*}
  S(5,5,5) &=
  \frac{7}{73728}\pi^4\zeta(11)+\frac{35}{24576}\pi^2\zeta(13)+\frac{63}{8192}\zeta(15),\\
  S(7,7,7) &=
  \frac{31}{35389440}\pi^6\zeta(15)+\frac{49}{1966080}\pi^4\zeta(17)+\frac{77}{262144}\pi^2\zeta(19)
  +\frac{429}{262144}\zeta(21).
\end{align*}

The values of $R(5,5,5)$, $R(7,7,7)$ and $R(9,9,9)$ listed
in~\cite{Tsu2} appear to be incorrect.  They should be
\begin{align*}
  R(5,5,5) &= \frac{16375}{147456}\pi^4\zeta(11) +
  \frac{573335}{49152}\pi^2\zeta(13)-\frac{2064195}{16384}\zeta(15),\\
  R(7,7,7) &= \frac{1048543}{70778880}\pi^6\zeta(15)+\frac{7339969}{3932160}\pi^4\zeta(17)
  +\frac{80740121}{524288}\pi^2\zeta(19)-\frac{899676921}{524288}\zeta(21),\\
  R(9,9,9) &=
  \frac{13421747}{7046430720}\pi^8\zeta(19)+\frac{738197141}{2113929216}\pi^6\zeta(21)
  +\frac{1919313253}{67108864}\pi^4\zeta(23)\\
  &+\frac{143948506845}{67108864}\pi^2\zeta(25)
  -\frac{1631416447375}{67108864}\zeta(27).
\end{align*}


\begin{thebibliography}{99} \raggedright


\bibitem{Berndt1} B.~Berndt, \textit{Ramanujan's Notebooks Part
I}, Springer, New York, 1985. [MR 0781125] (86c:01062)

\bibitem{ParmEuler} D.~Borwein, J.~M.~Borwein, and D.~M.~Bradley,
Parametric Euler sum identities, \textit{J.\ Math.\ Anal.\ Appl.},
\textbf{316} (2006), no.~1, 328--338.  doi:
10.1016/j.jmaa.2005.04.040 [MR 2201764] (2007b:11132) {\tt
http://arxiv.org/abs/math.CA/0505058}



\bibitem{GoldVar} J.~M.~Borwein and D.~M.~Bradley, Thirty-two
Goldbach variations, \textit{Internat.\ J.\ Number Theory},
\textbf{2} (2006), no.~1, 65--103.  doi: 10.1142/S1793042106000383
[MR 2217795] (2007e:11109)



\bibitem{BBB} J.~M.~Borwein, D.~M.~Bradley and D.~J.~Broadhurst,
{Evaluations of $k$-fold Euler/Zagier sums: a compendium of results
for arbitrary $k$},
%\textit{Electronic J.~Combinatorics},
\textit{Electron.~J.~Combin.}, \textbf{4} (1997), no.~2, \#R5. Wilf
Festschrift. [MR 1444152] (98b:11091)


\bibitem{BBBLc} J.~M.~Borwein, D.~M.~Bradley, D.~J.~Broadhurst and P.~Lison\v ek,
{Combinatorial aspects of multiple zeta values},
\textit{Electron.~J.~Combin.}, \textbf{5} (1998), no.~1, \#R38. [MR
1637378] (99g:11100) {\tt http://arXiv.org/abs/math.NT/9812020}

\bibitem{BBBLa} \bysame, %J.~M.~Borwein,  D.~M.~Bradley, D.~J.~Broadhurst and P.~Lison\v ek,
Special values of multiple polylogarithms, \textit{Trans.\ Amer.\
Math.\ Soc.},  \textbf{353} (2001), no.~3, 907--941. [MR 1709772]
(2003i:33003) {\tt http://arXiv.org/abs/math.CA/9910045}


\bibitem{BowBradSurvey} D.~Bowman and D.~M.~Bradley,
Multiple polylogarithms: a brief survey, \textit{Proceedings of a
Conference on $q$-Series with Applications to Combinatorics, Number
Theory and Physics}, (B.~C.~Berndt and K.~Ono eds., Urbana, IL,
2000), Contemporary Math., \textbf{291}, Amer.\ Math.\ Soc.,
Providence, RI, 2001, pp.\ 71--92.  [MR 1874522] (2003c:33021) {\tt
http://arXiv.org/abs/math.CA/0310062}

\bibitem{BowBrad3} \bysame, %D.~Bowman and D.~M.~Bradley,
The algebra and combinatorics of shuffles and multiple zeta values,
\textit{J.\ Combin.~Theory, Ser.\ A}, \textbf{97} (2002), no.~1,
43--61. [MR 1879045] (2003j:05010) {\tt
http://arXiv.org/abs/math.CO/0310082}

\bibitem{BowBrad1} \bysame, %D.~Bowman and D.~M.~Bradley,
Resolution of some open problems concerning multiple zeta
evaluations of arbitrary depth, \textit{Compositio Math.},
\textbf{139} (2003), no.~1, 85--100. doi:
10.1023/B:COMP:0000005036.52387.da [MR 2024966] (2005f:11196) {\tt
http://arXiv.org/abs/math.CA/0310061}


\bibitem{BowBradRyoo} D.~Bowman, D.~M.~Bradley, and J.~Ryoo,
{Some multi-set inclusions associated with shuffle convolutions and
multiple zeta values}, \textit{European J.~Combin.}, \textbf{24}
(2003), no.~1, 121--127. [MR 1957970] (2004c:11115)

\bibitem{Prtn} D.~M.~Bradley,
{Partition identities for the multiple zeta function}, \textit{Zeta
Functions, Topology, and Quantum Physics}, Developments in
Mathematics, \textbf{14}, T.\ Aoki, S.\ Kanemitsu, M.\ Nakahara, Y.\
Ohno (eds.) Springer-Verlag, New York, 2005,  pp.\ 19--29. ISBN:
0-387-24972-9. [MR 2179270] (2006f:11105) {\tt
http://arXiv.org/abs/math.CO/0402091}

\bibitem{DBqMzv} \bysame, %D.~M.~Bradley,
Multiple $q$-zeta values, \textit{J.\ Algebra}, \textbf{283} (2005),
no.~2, 752--798.  doi: 10.1016/j.jalgebra.2004.09.017 [MR 2111222]
(2006f:11106) {\tt http://arXiv.org/abs/math.QA/0402093}

\bibitem{DBqKarl} \bysame, %D.~M.~Bradley,
Duality for finite multiple harmonic $q$-series, \textit{Discrete
Math.}, \textbf{300} (2005), no.~1--3, 44--56. doi:
10.1016/j.disc.2005.06.008 [MR 2170113] (2006m:05019) {\tt
http://arXiv.org/abs/math.CO/0402092}



\bibitem{DBqSum} \bysame, %D.~M.~Bradley,
{On the sum formula for multiple $q$-zeta values}, \textit{Rocky
Mountain J.\ Math.}, to appear. {\tt
http://arxiv.org/abs/math.QA/0411274}

\bibitem{DBqDecomp} \bysame, {A $q$-analog of Euler's decomposition
formula for the double zeta function}, \textit{Internat.\ J.\ Math.\
Math.\ Sci.}, \textbf{2005} (2005), no.~21,  3453--3458.
doi:10.1155/IJMMS.2005.3453 [MR 2206867] (2006k:11174) {\tt
http://arxiv.org/abs/math.NT/0502002}

\bibitem{BK1}
D.~J.~Broadhurst and D.~Kreimer, Association of multiple zeta values
with positive knots via Feynman diagrams up to 9 loops,
\textit{Phys.\ Lett. B}, \textbf{393} (1997) no.~3-4, 403--412. [MR
1435933] (98g:11101)


\bibitem{Brown} F.~C.~Brown,
P\'eriodes des espaces des modules $\overline{\mathfrak{M}}\sb
{0,n}$ et valeurs z\^{e}tas multiples [Multiple zeta values and
periods of the moduli spaces $\overline{\mathfrak{M}}\sb {0,n}$]
\textit{C.\ R.\ Math.\ Acad.\ Sci.\ Paris}, \textbf{342} (2006),
no.~12, 949--954. [MR 2235616]

\bibitem{Cartier} P.~Cartier, Fonctions polylogarithmes, nombres
polyz\^{e}tas et groupes pro-unipotents, \textit{Ast\'erisque},
\textbf{282} (2002), viii, 137--173, S\'eminaire Bourbaki, 53'eme
ann\'ee, 2000--2001, Exp.\ No.~885. [MR 1975178] (2004i:19005)

\bibitem{Cartier06} \bysame, %P.~Cartier,
Values of the $\zeta$-function, \textit{Surveys in Modern
Mathematics}, 260--273, London Math.\ Soc.\ Lecture Note Ser.,
\textbf{321}, Cambridge Univ.\ Press, Cambridge, 2005. [MR 2166932]
(2006g:11182)


\bibitem{LE} L.~Euler, \textit{Meditationes Circa Singulare
Serierum Genus}, Novi Comm.\ Acad.\ Sci.\ Petropol., \textbf{20}
(1775), 140--186. Reprinted in ``Opera Omnia,'' ser.~I,
\textbf{15}, B.\ G.\ Teubner, Berlin (1927), 217--267.

\bibitem{LE2} \bysame, \textit{Briefwechsel}, vol.~1, Birh\"auser,
Basel, 1975. [MR 0497632]


\bibitem{Goldbach} L.~Euler and C.~Goldbach, \textit{Briefwechsel}
1729--1764, Akademie-Verlag, Berlin, 1965.

\bibitem{FlajSalv} P.\ Flajolet and B.\ Salvy, Euler sums and
contour integral representations, \textit{Experiment.\ Math.},
\textbf{7} (1998), no.~1, 15--35. [MR 1618286] (99c:11110)


\bibitem{Hoff02} M.~E.~Hoffman,
Periods of mirrors and multiple zeta values, \textit{Proc.\ Amer.\
Math.\ Soc.}, \textbf{130} (2002), no.~4, 971--974. [MR 1873769]
(2002k:14068)

\bibitem{HWZ} J.\ G.\ Huard, K.\ S.\ Williams and N.\ Y.\ Zhang, On
Tornheim's double series, Acta Arith., \textbf{75} (1996), no~2,
105--117. [MR 1379394] (97f:11073)


\bibitem{IT} K.~Ihara and T.~Takamuki, The quantum $\mathfrak{g}_2$
invariant and relations of multiple zeta values, \textit{J.\ Knot
Theory Ramifications}, \textbf{10} (2001), no.~7, 983--997.  [MR
1867104] (2002m:57016)

\bibitem{LeM} T.~Q.~T.~Le and J.~Murakami,
%Tu Quoc Thang Le and Jun Murakami,
Kontsevich's integral for the
Homfly polynomial and relations between values of multiple zeta
functions,
%\textit{Topology and its Applications},
\textit{Topology Appl.}, \textbf{62} (1995), no.~2, 193--206. [MR
1320252] (96c:57017)

\bibitem{SitSub} M.\ V.\ Subbarao and R.\ Sitaramachandra Rao, On
some infinite series of L.\ J.\ Mordell and their analogues,
\textit{Pacific J.\ Math.}, \textbf{119} (1985), no.~1, 245--255.
[MR 0797027] (87c:11091)

\bibitem{Tera1} T.~Terasoma, Selberg integrals and multiple zeta
values, \textit{Compositio Math.}, \textbf{133} (2002), no.~1,
1--24.  [MR 1918286] (2003k:11142)

\bibitem{Tera2} \bysame, %T.~Terasoma,
Mixed Tate motives and multiple zeta values, \textit{Invent.\
Math.}, \textbf{149} (2002), no.~2, 339--369. [MR 1918675]
(2003h:11073)

\bibitem{Tera3} \bysame, %T.~Terasoma,
Period integrals and multiple zeta values, \textit{S\={u}gaku},
\textbf{57} (2005), no.~3, 255--266.  [MR 2163671] (2006h:11110)


\bibitem{Torn} L.~Tornheim, Harmonic double series, Amer.\ J.\
Math., \textbf{72} (1950), 303--314. [MR 0034860] (11,654a)

\bibitem{Tsu1} H.~Tsumura, On some combinatorial relations for
Tornheim's double series, Acta Arith., \textbf{105} (2002), no.~3,
239--252. [MR 1931792] (2003i:11134)

\bibitem{Tsu2} \bysame, %H.~Tsumura,
On alternating analogues of Tornheim's double series, \textit{Proc.\
Amer.\ Math.\ Soc.}, \textbf{131} (2003), no.~12, 3633--3641. [MR
1998168] (2004e:11102)

\bibitem{Tsu3} \bysame, %H.~Tsumura,
Evaluation formulas for Tornheim's type of alternating double
series, Math.\ Comp., \textbf{73} (2004), no.~245, 251--258. [MR
2034120] (2005d:11137)



\bibitem{Wald02} M.~Waldschmidt,
Valeurs z\^{e}ta multiples: une introduction,
%\textit{Journal de Theorie des Nombres de Bordeaux},
\textit{J.~Th\'eor.~Nombres Bordeaux}, \textbf{12} (2000), no.~2,
581--595.  [MR 1823204] (2002a:11106)

\bibitem{Wald04} \bysame, %M.~Waldschmidt,
Multiple polylogarithms: an introduction, in \textit{Number Theory
and Discrete Mathematics}, (Chandigarh, 2000), Trends Math.,
Birkh\"auser, Basel, 2002, pp.\ 1--12.
%Hindustan Book Agency and Birkh\"auser Verlag, 2002, 1--12.
[MR 1952273] (2004d:33003)


\bibitem{Zud} V.~V.~Zudilin, Algebraic relations for
multiple zeta values (Russian), \textit{Uspekhi Mat.\ Nauk},
\textbf{58} (2003), no.~1, 3--32; translation in \textit{Russian
Math.\ Surveys}, \textbf{58} (2003), vol.~1, 1--29. [MR 1992130]
(2004k:11150)

\end{thebibliography}
\end{document}